\input amstex
\documentstyle{amsppt}
\input bull-ppt


\define\reals{\operatorname{\Bbb R}}

\define\ve{\varepsilon}

\define\KR{\text{Kr}}
\define\mr{\text{mr}}
\define\ov{\overline}

\topmatter
\cvol{26}
\cvolyear{1992}
\cmonth{Jan}
\cyear{1992}
\cvolno{1}
\cpgs{29-52}

\title Perspectives on Information-Based Complexity
\endtitle
\author J. F. Traub and H. Wo\'zniakowski \endauthor
\address \RM{(J. F. Traub)}
Department of Computer Science, Columbia University, New 
York,
New York 10027\endaddress
\address \RM{(H. Wo\'zniakowski)}
Department of Computer Science, Columbia University
and Institute of Informatics, University of Warsaw, 
Warsaw, Poland\endaddress
\date April, 1991\enddate
\endtopmatter

\document
\baselineskip=12pt
\TagsOnRight

\heading  1. \smc Introduction \endheading

Computational complexity studies the intrinsic difficulty 
of 
mathematically posed problems and seeks optimal means 
for their solutions.  This is a rich and diverse field; 
for the purpose of this paper we present a greatly 
simplified picture. 
\footnote""{This research was supported in part by the 
National Science
Foundation.}

Computational complexity may be divided into two 
branches, discrete and 
continuous. Discrete computational complexity studies 
problems such as graph
theoretic, routing, and discrete optimization; see, for 
example, 
Garey and Johnson [79]. Continuous computational 
complexity studies 
problems such as ordinary and 
partial differential equations, multivariate integration, 
matrix multiplication, and systems of polynomial 
equations. Discrete
computational complexity often uses the Turing machine 
model whereas 
continuous computational complexity tends to use the real 
number model. 
 
Continuous computational complexity may again be split 
into two branches.
The first deals with problems for which the information 
is {\it 
complete}. Problems where the information may be 
 complete are those 
for which the input is specified by a finite number of 
parameters. 
Examples include linear algebraic systems, matrix 
multiplication, and
systems of polynomial equations. Recently, Blum, Shub and 
Smale [89] 
obtained the first NP-completeness result over the reals 
for a problem
with complete information.

The other branch of continuous computational complexity 
is {\it
information-based complexity}, which is denoted for 
brevity as IBC.
Typically, IBC studies infinite-dimensional problems. 
These are
problems where either the input or the output are 
elements of
infinite-dimensional spaces.  Since digital computers can 
handle only
finite sets of numbers, infinite-dimensional objects such 
as functions
on the reals must be replaced by finite sets of numbers. 
Thus,
complete information is not available about such objects. 
 Only {\it
partial} information is available when solving an 
infinite-dimensional
problem on a digital computer.  Typically, information is 
{\it
contaminated} with errors such as round-off error, 
measurement error,
and human error.  Thus, the available information is 
partial and/or
contaminated.

We want to emphasize this point for it is central to IBC. 
{\it  Since only
partial and/or contaminated 
information is available, we can solve the original 
problem 
only approximately. The goal of IBC is to compute such an 
approximation 
as inexpensively as possible. }

In Figure 1 (see p.\ 30) we schematize the structure of 
computational complexity
described above. 

\fgh{9 pc}\caption{\smc Figure 1}

Research in the spirit of IBC was initiated in the Soviet 
Union by
Kolmogorov in the late 1940s. Nikolskij [50], then a 
graduate student
of Kolmogorov, studied optimal quadrature. 
This line of research was greatly advanced by Bakhvalov;
see, e.g., Bakhvalov [59, 64, 71].  
In the United States research in the spirit of IBC was 
initiated
by Sard [49] and Kiefer [53]. 
Kiefer reported the 
results of his 1948 MIT Master's Thesis that Fibonacci 
sampling is
optimal when approximating the maximum of a unimodal 
function. 
Sard studied optimal quadrature. 
Golomb and Weinberger [59] studied optimal approximation 
of linear 
functionals. Schoenberg [64] realized the close 
connection between
splines and algorithms optimal in the sense of Sard.

IBC is formulated as an abstract theory and it has 
applications in numerous
areas. The reader may consult TWW [88]
\footnote{When one of us is a coauthor, 
the citation will be made using only initials.} for some 
of the applications. 
IBC has benefitted from research in many fields. 
Influential have been 
questions, concepts, and results from complexity theory, 
algorithmic 
analysis, applied mathematics, numerical analysis, 
statistics, and the
theory of
approximation (particularly the work on $n$-widths and 
splines).

In this paper we discuss, in particular,
IBC research for two problems of numerical
analysis. We first 
contrast IBC and numerical analysis, limiting ourselves 
to just one
characteristic of each.

IBC is a branch of computational complexity, and optimal 
(or almost optimal) 
information and 
algorithms are obtained from the theory. In numerical 
analysis, 
particular classes of algorithms are carefully analyzed 
to see if they satisfy 
certain criteria such as convergence, error bounds, 
efficiency, and 
stability. 

Numerical analysis and IBC have different views on the 
problems which 
lie in their common domain. The authors of this paper 
have worked in
both numerical analysis and IBC, and believe the 
viewpoints are not 
right or wrong, just different. 

On the other hand, in many research groups around the 
world, people 
work on both numerical analysis and IBC, and do not draw 
a sharp distinction 
between the two. They believe IBC can serve as part of 
the theoretical 
foundation of numerical analysis. 

We believe there might be some profit in discussing the 
views of 
numerical analysis and IBC. Unfortunately Parlett [92]
\footnote{Citation to this paper will be made using only 
an initial.}
does not serve this purpose since, as we shall show, this 
paper
ignores relevant literature and is mistaken on issues of 
complexity
theory.

For example, P [92] contains a central misconception 
about IBC which
immediately invalidates large portions of the paper.  P 
[92] assumes
that the information is specified (or fixed).  Indeed, 
the first
``high level criticism'' is that IBC ``is not complexity 
theory'' (see
P [92, 2.A]), since ``specified information'' is used.

{\it But it is the very essence of IBC that both the 
information and the 
algorithms are varied. Indeed, one of the central 
problems of IBC is 
the optimal choice of information.} Significant portions 
of three 
monographs, TW [80] and TWW [83, 88], all of which are 
cited in P~[92], 
are devoted to this issue. We return to this issue in
\S 3 after notation has been established. 

In P [92], the author limits himself to ``matrix 
computations, which
is the area we understand best.''  We do not object to 
discussing
matrix computations, although they constitute a small 
fraction and are
atypical of IBC.  For example, in the recent monograph 
TWW [88], some
ten pages, just 2\%, are devoted to matrix computations. 
 Matrix
computations are atypical since complete information can 
be obtained
at finite cost. However, even in this particular area, P 
[92] ignores
relevant literature and does not exhibit a grasp of the 
complexity
issues. Since the discussion will, of necessity, assume 
some rather
technical details concerning matrix computations, we will 
defer it to
\S\S 5 and~6.

We stress that we are not questioning the importance of 
matrix
computations. On the contrary, they play a central role 
in scientific
computation. Furthermore, we believe there are some nice 
results and
deep open questions regarding matrix computations in IBC.

But the real issue is, after all, IBC in its entirety. 
P [92] is merely using the two
papers TW [84] and Kuczy\'nski [86] on matrix computations 
to criticize all of IBC. We therefore respond to general
criticisms in \S\S 3 and 4.  

To make this paper self-contained we 
briefly summarize the basic concepts of IBC in \S2. 
Section 7 deals with possible refinements of IBC.
A summary of our rebuttal to criticisms in P [92] is 
presented in \S8. 

\heading  2. Outline of IBC
\endheading
In this section we introduce the basic concepts of IBC 
and define the notation
which will be used for the remainder of this paper. We 
illustrate the
concepts with the example of multivariate integration, a 
typical application
of IBC. A more detailed account may be found in TWW [88]. 
Expository 
material may be found in W [85], PT [87], PW [87], and TW 
[91]. 
Let
$$S\: F\ \to \ G, 
$$
where $F$ is a subset of a linear space and $G$ is a 
normed linear
space. We wish to compute an approximation to $S(f\,)$ 
for all $f$ from $F$.

Typically, $f$ is an element from an infinite-dimensional 
space
and it cannot be represented on a digital computer. We 
therefore
assume that only {\it partial} information
\footnote{For simplicity, we will not consider 
contaminated  information in this paper.}
about $f$ is available. 
We gather this partial information about $f$ by computing 
information operations $L(f\,)$, where $L\in \Lambda$. Here
the class $\Lambda$ denotes a collection of information 
operations
that may be computed. We illustrate these concepts by an
example.

\dfn{Example: Multivariate integration} Let $F$ be a unit 
ball of the
Sobolev class $W_p^{r,d}$ of real functions defined on the
$d$-dimensional cube $D=[0,1]^d$ whose $r$th distributional
derivatives exist and are bounded in $L_p$ norm.  Let 
$G=\reals$ and
$$ S(f\,)=\int_Df(t)\,dt.  $$ Assume $pr\,>\,d$.  To 
approximate
$S(f\,)$, we assume we can compute only function values. 
 That is, the
class $\Lambda$ is a collection of $L\:F\to \reals$, such 
that for some
$x$ from $D$, $L(f\,)=f(x)$, $\forall f \in F$.\qed\enddfn

For each $f\in F$, we compute a number of information 
operations
from the class $\Lambda$. Let
$$
N(f\,)=[L_1(f\,),L_2(f\,),\dots,L_n(f\,)],\qquad L_i\in 
\Lambda,
$$
be the computed information about $f$. We stress that the 
$L_i$ as well 
as the number $n$ can be chosen adaptively. That is, the 
choice of $L_i$
may depend on the already computed 
$L_1(f\,),L_2(f\,),\dots,
L_{i-1}(f\,)$.
The number $n$ may also depend on the computed 
$L_i(f\,)$\<. (This
permits arbitrary termination criteria.)

$N(f\,)$ is called the information about $f$, and $N$ the 
information operator. 
In general, $N$ is many-to-one, and that is why 
it is impossible to recover the element~$f$, knowing 
$y=N(f\,)$ for $f\in F$.
 For this reason, the information $N$ is called {\it 
partial}. 

Having computed $N(f\,)$, we approximate $S(f\,)$ by an 
element
$U(f\,)=\phi(N(f\,))$, where $\phi\:N(F)\to G$. A mapping 
$\phi$ is called an
algorithm. 

The definition of error of the approximation $U$ depends 
on the setting. 
We restrict ourselves here to only two settings. In the 
worst case setting
$$
e(U)=\sup_{f\in F}\|S(f\,)-U(f\,)\|,
$$
and in the average case setting, given a probability 
measure $\mu$ on
$F$, 
$$
e(U)=\bigg(\int_F\|S(f\,)-U(f\,)\|^2\mu(df\,)\bigg)^{1/2}.
$$

\dfn\nofrills{Example \RM{(continued).}}
\ The information is given by
$$
N(f\,)=[f(x_1),f(x_2),\dots,f(x_n)]
$$
with the points $x_i$ and the number $n$ adaptively 
chosen. An example 
of an algorithm is a linear algorithm given by 
$U(f\,)=\phi(N(f\,))=\sum_{i=1}^na_i\,f(x_i)$
for some numbers $a_i$. 

In the worst case setting, the error is defined as the 
maximal distance
$|S(f\,)-U(f\,)|$ in the set $F$. In the average case 
setting,
the error is the $L_2$ mean of $|S(f\,)-U(f\,)|$ with 
respect to the 
probability measure $\mu$. The measure $\mu$ is sometimes
taken as a truncated Gaussian measure. \qed\enddfn

To define the computational complexity we need a model of 
computation. It is defined by two assumptions:
\roster
\item  We are charged for each information operation. 
That is, 
for every $L \in \Lambda$ and for every $f \in F$, the 
computation of $L(f\,)$ 
costs $c$, where $c$ is positive and fixed, independent 
of $L$ and $f$. 
\item  Let $\Omega$ denote the set of permissible 
combinatory
operations including the addition of two elements in $G$, 
multiplication by a scalar in $G$, arithmetic operations, 
comparison
of real numbers, and evaluations of certain elementary 
functions.
We assume that each combinatory operation is performed 
exactly 
with unit cost. 
\endroster

In particular, this means that we use the real number 
model, where 
we can perform operations 
on real numbers exactly and at unit cost. Modulo roundoffs 
and the very important concept of numerical stability, 
this corresponds to
floating point arithmetic widely used for solving 
scientific 
computational problems. 

We now define the cost of the approximations. Let 
$\text{cost}(N,f\,)$ denote
the cost of computing the information $N(f\,)$. Note that 
$\text{cost}(N,f\,)
\ge c\,n$, and the inequality may occur since adaptive 
selection of
$L_i$ and $n$ may require some combinatory operations. 

Knowing $y=N(f\,)$, we compute $U(f\,)=\phi(y)$ by 
combining the information
$L_i(f\,)$. Let $\text{cost}(\phi,y)$ denote the
number of combinatory operations from $\Omega$ needed to 
compute $\phi(y)$. 
We stress that $\text{cost}(N,f\,)$ or 
$\text{cost}(\phi,y)$ may be equal
to infinity if $N(f\,)$ or $\phi(y)$ use an operation 
outside $\Omega$
or infinitely many operations from $\Lambda$ or $\Omega$, 
respectively. 

The cost of computing $U(f\,)$, cost$(U,f\,)$, is given by
$$
\text{cost}(U,f\,)\,=\,\text{cost}(N,f\,)+
\text{cost}(\phi,N(f\,)).
$$
Depending on the setting, the cost of $U$ is defined as 
follows.
In the worst case setting
$$
\text{cost}(U)\,=\,\sup_{f \in F}\,\text{cost}(U,f\,),
$$
and in the average case setting
$$
\text{cost}(U)\,=\,\int_F\text{cost}(U,f\,)\,\mu(df\,).
$$

We are ready to define the basic notion of 
$\ve$-complexity.
The $\ve$-complexity is defined as the minimal cost among 
{\it all}
$U$ with error at most $\ve$,
$$
\text{comp}(\ve)\,=\,\inf\{\text{cost}(U)\: U\ \text{such 
that}\ 
e(U)\le \ve\}.
$$
(Here we use the convention that the infimum of the empty 
set is taken to be
infinity.) Depending on the setting, this defines
the worst case or average case $\ve$-complexity.

We stress that we take the infimum over {\it all} 
possible $U$ for which 
the error does not exceed $\ve$. Since $U$ can be 
identified with the pair
$(N,\phi)$, where $N$ is the information and $\phi$ is the 
algorithm that uses that information, this means that we 
take the 
infimum over {\it all} information $N$ consisting of 
information operations
from the class $\Lambda$, and over {\it all} algorithms 
$\phi$ that use
$N$ such that $(N,\phi)$ computes approximations with 
error at most $\ve$. 
 
\rem{Remark}
The complexity depends on the set $\Lambda$ of permissible 
information operations and on the set $\Omega$ of 
permissible 
combinatory operations. Both sets are necessary to define 
the 
complexity of a problem. This is beneficial because 
the dependence of complexity on $\Lambda$ and 
$\Omega$ enriches the theory; it enables us to study the 
power of 
specified information or combinatory operations. 
We illustrate the role of $\Lambda$ and $\Omega$ by a 
number of examples. 

We begin with the role of $\Lambda$.
Assume that $F$ is a subset of a linear space of functions.
Let $\Lambda_1$ consist of all linear functionals, and 
let $\Lambda_2$ 
consist of function evaluations. For many applications 
$\Lambda_2$ is 
more practical. Let $\Omega$ be defined as above.

Consider the integration example. For this problem, 
$\Lambda_1$ is not
a reasonable choice 
since any integral could be computed exactly with cost 
$c$. 
For $\Lambda_2$, we get the multivariate integration 
problem  
discussed in this section. 

Consider next the approximate solution of $2m$\<th-order 
elliptic linear partial 
differential equations whose right-hand side belongs to 
the unit ball
of $H^{r}(D)$ for a bounded simply-connected $C^{\infty}$ 
region $D$ of 
$\reals^d$. Let $G=H^m(D)$. Werschulz has shown that the 
worst case complexity 
in the class $\Lambda_1$ is proportional to $\ve^{-d/(r+
m)}$, and in the 
class $\Lambda_2$ it is proportional to $\ve^{-d/r}$; 
a thorough study of this subject may be found  in the 
research 
monograph Werschulz [91]. 
Thus, the complexity
penalty for using $\Lambda_2$ rather than $\Lambda_1$ 
goes to infinity as 
$\ve$ goes to zero for $m>0$; see also TWW [88, Chapter 
5, Theorem 5.9]. 
On the other hand, Werschulz has shown that 
the complexity of Fredholm integral equations of the 
second  kind is roughly the same for $\Lambda_1$ and 
$\Lambda_2$; see Werschulz [91] as well as TWW [88, 
Chapter 5, \S 6]. 

We now illustrate the role of $\Omega$ for the 
approximate solution of
scalar complex polynomial equations of degree $d$ using 
complete
information, i.e., $\Lambda$ consists of the identity 
mapping.  Let
$\Omega_1$ consist of the four arithmetic operations 
(over the complex
field), and let $\Omega_2$ consist of the four arithmetic 
operations
and complex conjugation. We confine ourselves to purely 
iterative
algorithms.  Then for $d\ge 4$, McMullen [85] proved that 
the problem
{\it cannot} be solved for the class $\Omega_1$, whereas 
Shub and
Smale [86] proved that the problem {\it can} be solved 
for the class
$\Omega_2$.  The positive result of Shub and Smale [86] 
also holds for
systems of complex multivariate polynomials of degree 
$d$. Hence, the
arithmetic operations are too weak for approximate 
polynomial zero
finding, whereas also permitting complex conjugation 
supplies enough
power to solve the problem. \qed\endrem

\dfn\nofrills{Example \RM{(continued).}}
\ For the integration problem, 
the model of computation states that one function 
evaluation costs $c$, and each arithmetic operation, 
comparisons of real
numbers, and evaluations of certain elementary functions
can be performed exactly at unit cost. Usually $c\gg 1$. 

The worst case $\ve$-complexity for the unit ball of 
$W_p^{r,d}$ is 
as follows. 
For $pr\,>\,d$,
$$
\text{comp}(\ve)=\Theta(c\ve^{-d/r})\quad \text{as}\ \ve 
\to 0;
$$
see Novak [88] for a recent survey. Take $p=+\infty$. 
Then for $d$
large relative to $r$, the worst case $\ve$-complexity is 
huge even
for moderate $\ve$. Furthermore, if only continuity of 
functions is
assumed, then the problem cannot be solved since
$\text{comp}(\ve)=+\infty$.

For the average case setting, let $F$ be 
the unit ball in the $\sup$ norm of continuous functions. 
Let $\mu$ be 
a truncated classical Wiener sheet measure; see, e.g., 
TWW [88, p.\ 218].
Then using results from number theory concerning 
discrepancy
(see Roth [54, 80]),  we have
$$
\text{comp}(\ve)=\Theta(c\ve^{-1}(\log\ve^{-1})
^{(d-1)/2}) \quad \text{as}\ \ve \to 0;
$$
see W [87, 91]. 
Thus, the average case complexity depends only mildly on 
the dimension~$d$. 
(The same $\Theta$ result holds 
if the unit ball is replaced by the entire space of 
continuous functions.)
To get an approximation with cost proportional to 
$\text{comp}(\ve)$, 
it is enough to compute the arithmetic mean 
$n^{-1}\sum_{i=1}^nf(x_i)$,
where $n=\Theta(\ve^{-1}(\log \ve^{-1})^{(d-1)/2})$, and 
the points $x_i$
are derived from Hammersley points. \qed \enddfn

A goal of IBC is to find or estimate the 
$\ve$-complexity, and to find 
an  $\ve$-complexity optimal $U$, or equivalently, an
$\ve$-complexity optimal pair $(N,\phi)$.
By {\it $\ve$-complexity optimality} of $U$ we mean that 
the error of $U$ is 
at most $\ve$ and the cost of $U$ is equal to, or not 
much greater than, the 
$\ve$-complexity. For a number of problems this goal has 
been achieved due 
to the work of many researchers. 

Many computational problems can be formulated using the 
approach outlined 
above. For some problems, including the two matrix 
computation problems
discussed in P [92], we need a more general formulation. 
We now
briefly discuss this more general formulation; details 
can be found in  
TWW [83, 88]. 

Let $F$ and $G$ be given sets, and $W$ be a given mapping
$$
W\:F\times[0,+\infty)\ \to\ 2^G.
$$
We assume that $W(f,0)$ is nonempty and grows as $\ve$ 
increases, i.e.,
for any $\ve_1\le \ve_2$ we have $W(f,\ve_1)\subset 
W(f,\ve_2),\ 
\forall f \in F.$ 

We now wish to compute an element $U(f\,)$ which belongs 
to $W(f,\ve)$ 
for all $f \in F$. 
The definitions of  $U$ as well as the cost of $U$ 
are unchanged. The error of $U$ is now defined as 
follows. The 
error of $U$ for $f$ from $F$ is
$$
e(U,f\,)=\inf\{\eta:\ U(f\,)\in W(f,\eta)\}.
$$
Then the error of $U$ is defined as 
$e(U)=\sup_{f\in F}e(U,f\,)$ in the worst case setting, 
and 
$e(U)=(\int_Fe^2(U,f\,)\,\mu(df\,)\,)^{1/2}$ in the 
average case setting. 
Note that for
$$
W(f,\ve)=\{g\in G:\ \|S(f\,)-g\|\le \ve\}
$$
we have $e(U,f\,)=\|S(f\,)-U(f\,)\|$ and the two 
formulations coincide. 

Finally, we illustrate how the two matrix computation 
problems fit
in this formulation.

\subheading{\RM{(i)}\quad Large linear systems}
We wish to approximate the solution
of a large linear system $Az\!=\!b$ by computing a vector 
$x$ with residual 
at most $\ve\!$, $\|Ax-b\|\le \ve$. Here, $b$ is a given 
vector, $\|b\|=1$,
and $A$ belongs to a class $F$ of $n\times n$ nonsingular 
matrices. 
The vectors $x$ are computed by using matrix-vector 
multiplications $Az$ for
any vector $z$. 

This problem corresponds to taking $G=\reals^n$ and 
$$
W(A,\ve)=\{x\in G:\ \|Ax-b\|\le \ve\},\quad \forall A\in F.
$$
The class $\Lambda$ of information operations is now 
given by
$$\align
\Lambda=\{L:F\to \reals^n\: &\text{there exists a vector 
$z\in\reals^n$}\\
&\text{such that $L(A)=Az,\ \forall A\in F$}\}.
\endalign
$$

\subheading{\RM{(ii)}\quad Eigenvalue problem}

For a matrix $A$ from a class
$F$ of $n\times n$ symmetric 
matrices, we wish to compute an approximate eigenpair 
$(x,\lambda)$, where $x\in \reals^n$ with $\|x\|=1$, 
and $\lambda\in \reals$, such that
$$
\|Ax-\lambda\,x\|\ \le \ \ve\,\|A\|.
$$
As in (i), the pairs $(x,\lambda)$ are computed by using 
matrix-vector multiplications.

This problem corresponds to taking $G=B^n\times \reals$, 
where $B^n$ is the 
unit sphere of $\reals^n$, and 
$$
W(A,\ve)=\{(x,\lambda)\in G:\ \|Ax-\lambda\,x\|\le 
\ve\,\|A\|\},
\quad\forall A\in F.
$$
The class $\Lambda$ is the same as in (i).

\heading  3. The role of information
\endheading

Information is central to IBC. We indicate briefly why 
the distinction 
between information and algorithm is so powerful. We then 
respond to two 
general criticisms in P [92] regarding information. 

As explained in \S 2, the approximation $U(f\,)$ is 
computed by
combining information operations from the class 
$\Lambda$.  Let
$y=N(f\,)$ denote this computed information. In general, 
the operator
$N$ is many-to-one, and therefore the set $N^{-1}(y)$ 
consists of many
elements of $F$ that cannot be distinguished from $f$ 
using $N$.
Then the set $SN^{-1}(y)$ consists of all elements from 
$G$ which are
indistinguishable from $S(f\,)$. Since $U(f\,)$ is the 
same for any $f$
from the set $N^{-1}(y)$, the element $U(f\,)$ must serve 
as an
approximation to any element $g$ from the set 
$SN^{-1}(y)$. It is
clear that the quality of the approximation $U(f\,)$ 
depends on the
``size'' of the set $SN^{-1}(y)$. In the worst case 
setting, define
the {\it radius of information} $r(N)$ as the maximal 
radius of the
set $SN^{-1}(y)$ for $y\in N(F)$. (The radius of the set 
$A$ is the
radius of the smallest ball which contains the set $A$.)

Clearly, the radius of information $r(N)$ is a sharp 
lower bound on the 
worst case error of any $U$. We can guarantee an 
$\ve$-approximation
iff $r(N)$ does not exceed $\ve$ (modulo a technical 
assumption 
that the corresponding infimum is attained). 

The cost of computing $N(f\,)$ is at least $cn$, where 
$c$ stands for the cost
of one information operation, and $n$ denotes their 
number in 
the information $N$. By the {\it $\ve$-cardinality number 
$m(\ve)$}
we mean the minimal number $n$ 
of information operations for which the information $N$ 
has 
radius $r(N)$ at most equal to $\ve$. 
From this we get a lower bound on the $\ve$-complexity in 
the worst 
case setting,
$$
\text{comp}(\ve)\,\ge\,cm(\ve).
$$
For some problems (see TWW [88, Chapter 5, \S 5.8]) 
it turns out that it is possible to find an information
operator $N_{\ve}$ consisting of $m(\ve)$ information 
operations, 
and a mapping $\phi_{\ve}$ such that the approximation
$U(f\,)=\phi_{\ve}(N_{\ve}(f\,))$ has error at 
most $\ve$ and $U(f\,)$ can be computed with cost at most 
$(c+2)\,m(\ve)$. This yields an upper bound on the 
$\ve$-complexity,
$$
\text{comp}(\ve)\,\le\,(c+2)\,m(\ve).
$$
Since usually $c\gg 1$, the last two inequalities yield 
the almost exact 
value of the $\ve$-complexity,
$$
\text{comp}(\ve)\,\simeq\,cm(\ve).
$$
This also shows that the pair $(N_{\ve},\phi_{\ve})$ is 
almost 
$\ve$-complexity optimal.

In each setting of IBC one can define a radius of 
information such that
we can guarantee an $\ve$-approximation iff $r(N)$ does 
not exceed 
$\ve$; see TWW [88]. This permits one to obtain 
complexity bounds in 
other settings. 

What is the essence of this approach? The point is that 
the radius of 
information as well as the $\ve$-cardinality number 
$m(\ve)$ and 
the information $N_{\ve}$ 
do not depend on particular algorithms, 
and they can often be expressed 
entirely in terms of well-known mathematical concepts. 
Depending on the 
setting and on the particular problem, the radii of 
information,
the $\ve$-cardinality numbers, and the information 
$N_\ve$ are related 
to Kolmogorov and Gelfand $n$-widths, $\ve$-entropy, 
the traces of correlation operators of
conditional measures, discrepancy theory, the minimal 
norm of splines, etc. 

In summary, there are two reasons why one can sometimes 
obtain sharp bounds on 
$\ve$-complexity in IBC.  The first is the distinction 
between information 
and algorithm. The second is that, due to this 
distinction, 
one can draw on powerful results in pure and applied 
mathematics.

We now respond to two central criticisms in P [92] 
regarding information. 
He asserts:
\roster 
\item"{(i)}" The information is specified (or given) and 
therefore this ``is 
not complexity theory;'' see P [92, 2.A].
\item"{(ii)}" There is an ``artificial distinction 
between information and 
algorithm;''  see P [92, 1].
\endroster

(i) P [92] repeatedly asserts that the information 
is ``specified'' or ``given.'' We have already referred 
to this misconception 
in our introduction and will amplify our response here. 

Varying the information {\it and} the algorithms is 
characteristic of IBC.
(For problems for which information is complete,
i.e., $N$ is one-to-one, only the algorithms can be 
varied.) The definition
of computational complexity in our work always entails 
varying both 
information and algorithms; see, for example, TW [80, 
Chapter 1, 
Definition 3.2], TWW [83, Chapter 5, \S 3], W [85, 2.5], 
PW [87, II], TWW [88, Chapter 3, \S 3]. 

Furthermore the study of optimal information, which of 
course makes sense only
if the information is being varied, is a constant theme 
in our work; see,
for example, TW [80, Chapters 2 and 7], TWW [83, Chapter 
4], W [85, 3.5],
PW [87, III D,  V C], TWW [88, Chapter 4, \S 5.3, Chapter 
6, 
\S 5.5].

Here, we have responded to criticism (i) in general. In 
\S\S 5 and 6 
we respond for the case of matrix computations. 

(ii.1) P [92, 1] claims there is an ``artificial 
distinction between
information and algorithm.'' That is, he argues that 
writing the
approximation $U(f\,)=\phi(N(f\,))$ is sometimes 
restrictive. We are
surprised that he does not produce a single example to 
back his claim.

(ii.2) P [92, Abstract] states that 
``a sharp distinction is made between information and 
algorithms restricted 
to this information. Yet the information itself usually 
comes from an 
algorithm and so the distinction clouds the issues and 
can lead to true 
but misleading inferences.'' 

We once again explain our view of the issues involved 
here using a
simple integration example.

As in \S 2 assume that we can compute function values. 
How can we
approximate the integral of $f$? The approximation 
$U(f\,)$ can be computed 
by evaluating $f$ at a number of points, say at 
$x_1,x_2,\dots, x_n$, 
and then the computed values $f(x_i)$ are combined to get 
$U(f\,)$. 
Computations involving $f(x_i)$, the adaptive selection 
of the points $x_i$,
and the adaptive choice of $n$ 
constitute the information $N(f\,)$. 
Denoting by $\phi$ the mapping which combines $N(f\,)$, 
we get 
$U(f\,)=\phi(N(f\,))$. 

We do not understand why this is restrictive, why it 
clouds the
issues, and why it leads to ``true but misleading 
inferences.'' As explained
in the first part of this section,  the distinction 
between information 
and algorithm sometimes enables us to find sharp bounds 
on complexity.

\heading  4. The domain $F$
\endheading
A basic concept in IBC is the domain $F$. A central 
criticism of IBC 
in P [92] concerns~$F$. The assertion is that there are 
two 
difficulties with $F$:
\roster
\item"{(i)}"  There is no need for $F$.
\item"{(ii)}"  There should be a charge for knowing 
membership in $F$.
\endroster

Concerning (i), the second ``high level criticism'' P 
[92, 2.B]
states:

``The ingredient of IBCT that allows it to generate 
irrelevant results is the problem class~$F$. $F$ does not 
appear 
in our brief description of the theory in the second 
paragraph of 
\S1 because it is not a logically essential ingredient 
but rather 
a parameter within IBCT.'' 

Concerning (ii), P [92, Abstract] states:

``By overlooking $F$'s 
membership fee the theory sometimes distorts the 
economics of problem 
solving in a way reminiscent of agricultural subsidies.'' 

First, why is $F$ needed?

(i.1) The set $F$ is necessary since it is the domain of 
the operator $S$,
or part of the domain of the operator $W$.

One need not say anything further; an operator must have 
a domain.
Nevertheless we will add a few additional points 
regarding the domain $F$. 

(i.2) For discrete or finite-dimensional problems one can 
sometimes
take the ``maximal'' set as $F$. 
Thus, in studying the complexity of matrix multiplications 
one usually takes $F$ as the set of all $n\times n$ 
matrices. In 
graph-theoretic complexity one often takes $F$ as the set 
of all 
graphs $(V,E)$, where $V$ is the set of vertices and $E$ 
is the set of edges.

However, for infinite-dimensional problems one cannot 
obtain meaningful
complexity results if $F$ is too large. For example, the 
largest $F$ one might
take for integration is the set of Lebesgue-integrable 
 functions, but then
$\roman{comp}(\ve)=+\infty,\ \forall\,\ve\ge 0$ in the 
worst  case setting. The
$\ve$-complexity remains infinite even if $F$ is  the set 
of continuous
functions. 

To make the complexity of an infinite-dimensional problem 
finite, 
one must take a smaller $F$ in the worst case setting or 
switch 
to the average case setting. Thus, as we saw in \S2, in 
the average 
case setting with a Wiener measure, the complexity is 
finite even if $F$ is 
the set of continuous functions. 

(i.3) The use of $F$ is not confined to IBC. 
In discrete computational complexity researchers often 
use a set $F$ 
which is smaller than the maximal set. For example, if 
$F$ is the set of all
graphs then many problems are NP-complete. If $F$ is a 
specified smaller 
set, then depending on the problem it may remain 
NP-complete or it may be
solvable in polynomial time. See, for example, Garey and 
Johnson [79]. 

(i.4) We believe the dependence of complexity on $F$ is 
part of the 
richness of IBC. For example, in the integration problem 
it is interesting
to\ know how complexity depends on the number of 
variables and
the smoothness of the integrands. 

(i.5) For a moment, we specialize our remarks to matrix 
computations. One 
could study the complexity of large linear systems for 
the set $F$
of all invertible matrices of order $n$. Then to compute 
an $\ve$-approximation one would have to recover the 
matrix $A$ 
by computing $n$ matrix-vector multiplications; this is 
a negative result. 

We find  criticism (i) particularly odd since 
an entire book, Parlett [80], 
is devoted to only the eigenvalue problem for symmetric 
matrices.
The reason is, of course, that 
the algorithms and the analysis for the symmetric 
eigenvalue
problem are very different than for arbitrary matrices. 
But then
why is the concept of $F$ so elusive?

Researchers in numerical linear algebra often consider 
other 
important subsets of 
matrices such as tridiagonal, Toeplitz, or Hessenberg 
matrices.

We turn to the criticism that there should be a charge for 
knowing membership in $F$. 

(ii.1) Is IBC being held to a  higher standard?
Do researchers in other disciplines charge for $F$? For 
example, researchers 
in numerical analysis often analyze the cost and error of 
important 
algorithms. The analysis depends on $F$. To give a simple 
example,
the analysis of the composite trapezoidal rule usually 
requires that the second
derivative of the integrand is bounded. There is no 
charge for membership in $F$. Indeed, how would one 
charge for knowing 
that a function has a bounded second derivative?

(ii.2) We believe that P [92] confuses two different 
problems:
\roster
\item "{(a)}" approximation of $S(f\,)$ for $f$ from $F$,
\item "{(b)}" the domain membership problem; that is, 
does $f$ belong to $F$?
\endroster
Domain membership is an interesting problem
which may be formulated within the IBC framework, 
although it has nothing
to do with the original problem of approximating $S(f\,)$ 
for $f\in F$. 

We outline how this may be done. 
First, to make the domain membership problem meaningful 
we {\it must} 
define the domain of $f$, say the set $\ov F$, in such a 
way that 
the logical values of $f\in F$ vary with $f$ from $\ov F$,
i.e., $\emptyset \ne F\cap\ov F\ne \ov F$.  Let
$\ov S:\ov F\to \{0,1\}\subset \reals $ be given by
$$\ov S(f\,)=\chi_F(f\,),\quad \forall f \in \ov F,$$
where $\chi_F$ is the characteristic (indicator) function 
of $F$. 

Then the problem is to compute $\ov S(f\,)$ exactly or 
approximately.
Observe that we now {\it assume} that $f\in \ov F$ just 
as we 
{\it assumed} that $f\in F$ for problems of type (a). 

For the domain membership problem we charge for computing 
an 
approximation to $\ov S(f\,)$, and 
the complexity of the domain membership problem is the 
minimal cost
of verifying whether $f\in F$. 

In the worst case setting, only the exact computation of 
$\ov S(f\,)$
makes sense since 
for $\ve\ge \tfrac 12$ the problem is trivial, and for 
$\ve <\tfrac 12$ 
it is 
the same as for $\ve=0$. However for the average case or 
probabilistic 
settings, an $\ve$-approximation may be reasonable. 
For instance we may wish to compute  $\ov S(f\,)$ with 
probability
$1-\ve$. 

It is easy to see that, in general, the domain membership 
problem cannot be
solved in the worst case setting. To illustrate this, let 
$\ov F$
be the set of continous functions, and let~$F$ be
the set of $r$ times continuously
differentiable functions, $r\ge 1$. Let the class 
$\Lambda$ of 
information operations 
consist of function values. It is obvious that knowing $n$ 
values of~$f$, no matter 
how large $n$ may be, there is no way to verify whether
$f$ is a member of $F$. 

The domain membership problem can be studied in the 
average case or 
probabilistic settings. Its complexity may be large or 
small depending 
on $\ov F$ and $F$. 
An example of work for this problem is
Gao and Wasilkowski [90] who study a particular domain 
membership problem.

(ii.3)  
Finally, we are at a loss to understand the following 
sentence from
P [92, 2.B], ``Whenever $F$ is very large
(for example, the class of continuous functions or the 
class of invertible
matrices) then it is realistic to assign no cost to it.'' 
Why is it realistic 
to assign no cost for ``large'' $F$, and why is it 
necessary to assign 
cost to ``small'' $F$? Where is the  
magic line which separates large $F$ from small $F$?

\heading 5. Large linear systems
\endheading
We briefly describe IBC research on large linear systems 
and then
respond to the criticisms in P [92]. Let
$$Ax=b,$$
where $A\in F$, and $F$ is a class of $n\times n$ 
nonsingular matrices.
Here $b$ is a known $n\times 1$ vector normalized such 
that $\|b\|=1$, and 
$\|\cdot\|$ stands for the spectral norm. 

Our problem is defined as follows. For any $A\in F$ and 
any $\|b\|=1$
compute an $\ve$-approximation $x$,
$$\|Ax-b\|\le \ve.$$

Usually $A$ is sparse and therefore $Az$ can be computed 
in time and storage
proportional to $n$. 
It is therefore reasonable for large linear systems to 
assume 
that the class $\Lambda$ of information operations consists
of matrix-vector multiplications. That is, we can compute 
$Az_1,Az_2,\dots,
Az_k$, where $z_i$ may depend on the known vector $b$ and 
on the previously
computed vectors $Az_1,\dots,Az_{i-1}$. To stress that 
the right-hand 
side vector
$b$ is known we slightly abuse the notation of \S 2 and 
denote
$$
N_k(A,b)=[b,Az_1,\dots,Az_k],\qquad A\in F, \tag 5.1
$$
as the information about the problem. The number $k$ is 
called the 
{\it cardinality} of information. For this to be of 
interest, we need
$k\ll n$. 

{\it Krylov information} is the special case when we take 
$z_1=b$ and 
$z_i=Az_{i-1}$. Thus Krylov information is given by
$$
N_k^{\KR}(A,b)=[b,Ab,\dots,A^kb].
$$

In what follows we will use the concept of {\it 
orthogonal invariance} 
of the class $F$. The class $F$ is 
{\it orthogonally invariant} iff 
$$
A\in F\quad  \text{implies}\quad  Q^{\roman T}A\,Q\in F
$$
for any orthogonal matrix $Q$, i.e., satisfying 
$Q^{\roman T}Q=I$.

Examples of orthogonally invariant
classes include many of practical interest such as
symmetric matrices, symmetric positive definite matrices, 
and matrices with 
uniformly bounded condition numbers.

We first discuss {\it optimal} information for large 
linear 
systems which is defined
as follows. 
The $\ve$-cardinality number $m(\ve)$ (see \S3)
denotes now the minimal 
cardinality $k$ of {\it all} information $N_k$ of the 
form (5.1)
with $r(N_k)\le \ve$. Obviously, 
$m(\ve)$ depends on the class $F$ and the class $\Lambda$. 
The information $N_k^*$ is {\it optimal} 
iff $k=m(\ve)$ and $r(N_k^*)\le\ve$. 	

\rem{Remark}
In \S2 we define the $\ve$-complexity optimality of a 
pair $(N,\phi)$. 
In this section optimality of information $N^*_k$ is 
introduced.
How are these two optimality notions related? 

In general, they are not. However, as already indicated 
in \S2, for many
problems the cost of computing $N^*_k(A,b)$ is 
proportional 
to $cm(\ve)$ and  
there exists an algorithm $\phi^*$ that uses $N^*_k$ and 
has error $\ve$ and 
combinatory cost proportional to $m(\ve)$. Then the pair 
$(N^*_k,\phi^*)$ is (almost) $\ve$-complexity optimal. In 
this case,
the two notions of optimality coincide and 
the complexity analysis reduces to the problem of finding 
optimal information. Details may be found in TWW [88, 
Chapter 4,
\S4]. \qed\endrem

\smallskip

In  TW [84] we conjecture that for the class $\Lambda$ of 
matrix-vector multiplications and for any orthogonally 
invariant $F$, 
Krylov information is {\it optimal}. 

Chou [87], based on Nemirovsky and 
Yudin [83], shows that Krylov information is  {\it 
optimal} modulo
a multiplicative factor of 2. More precisely, 
let $m^{\KR}(\ve)$ denote the minimal cardinality $k$ of 
Krylov information for which $r(N^{\KR}_k)\le \ve$. For 
{\it any} orthogonally invariant class $F$, we have 
$$m(\ve)\le m^{\KR}(\ve)\le 2\,m(\ve)+2.
$$
Recently, Nemirovsky [91] shows that for a number of 
important 
orthogonally invariant classes $F$ and for
$m(\ve)\le \tfrac 12(n-3)$, Krylov information is {\it 
optimal},
$$
m(\ve)=m^{\KR}(\ve).
$$

We now discuss algorithms that use Krylov information. 
We recall the definition of
the classical minimal residual (mr) algorithm; see, e.g., 
Stiefel [58].
The mr algorithm, $\phi^{\mr}$, 
uses Krylov information $N_k^{\KR}(A,b)$ and computes 
the vector $x_k$ such that 
$$\align
\|Ax_k-b\|= \min\{\|W_k(A)b\|\: 
 &W_k \ \text{is a polynomial}\\ 
& \text{of degree} \le k \ \ \text{and}\ \  W_k(0)=1\}.
\endalign
$$
Thus, by {\it definition} the mr algorithm minimizes the 
residual
in the class of {\it polynomial} algorithms. 

The mr algorithm has many good properties. 
Let $m^{\KR}(\ve,\phi^{\mr})$ 
denote the minimal cardinality of 
Krylov information needed to compute an 
$\ve$-approximation by the mr 
algorithm. 
Obviously, $m^{\KR}(\ve)$ denotes the minimal cardinality
of Krylov information needed to compute an 
$\ve$-approximation in the class
of {\it all} algorithms.
For any orthogonally invariant class~$F$, we have (see TW 
[84])
$$
m^{\KR}(\ve)\le m^{\KR}(\ve,\phi^{\mr})\le m^{\KR}(\ve)+1.
$$
These bounds are sharp. That is, for some $F$ we have
$m^{\KR}(\ve)=m^{\KR}(\ve,\phi^{\mr})$, and for other~$F$ 
we have 
$m^{\KR}(\ve,\phi^{\mr})=m^{\KR}(\ve)+1$. 

For all practically important cases, $m^{\KR}(\ve)$ is 
large and there 
is no significant difference between 
$m^{\KR}(\ve,\phi^{\mr})$ and 
$m^{\KR}(\ve)$. Therefore the mr algorithm is {\it always} 
recommended as long as $F$ is orthogonally invariant. 

The mr algorithm minimizes, up to an additive term of 
$1$, the 
number of matrix-vector multiplications needed to compute 
an 
$\ve$-approximation among {\it all} algorithms that use 
Krylov
information in an orthogonally invariant class $F$. In 
this sense,
the mr algorithm is {\it Krylov-optimal}, or for brevity, 
{\it optimal}.

We comment on the mr algorithm.
\roster
\item "{(1)}" The mr algorithm  computes $x_k$ {\it 
without} using 
the additional properties of $A$, $A\in F$, 
given in the definition of the class $F$. This is 
desirable since the computation of $x_k$ is the same for 
all $F$.
The vector $x_k$ can be computed by the well-known 
three-term recurrence 
formula using at most $10\,kn$ arithmetic operations. 
\item "{(2)}" Although the mr algorithm
competes with {\it all} algorithms, in particular with 
algorithms
that may use the additional properties of $A$ 
given in the definition of $F$, the mr algorithm
can lose at most one insignificant step. Equivalently, 
one may say that 
{\it for any orthogonally invariant class $F$, 
the a priori information about the class $F$ and the fact 
that $A \in F$
is worth at most one step. }
\item "{(3)}" 
On the other hand, if $F$ is {\it not } orthogonally 
invariant then
the mr algorithm may lose its good properties. Example 
 3.5 
of TW [84] provides such a class for which the worst 
happens;
the mr algorithm takes $n$ steps to solve the problem, 
whereas
the optimal algorithm, which is nonpolynomial, takes only 
one step. 
\endroster

For an orthogonally invariant class $F$ and for the class 
$\Lambda$ of
matrix-vector multiplications, these results yield that 
the pair Krylov information and mr algorithm is (almost) 
$\ve$-complexity optimal in the sense of \S 2. 
Furthermore, we 
have rather tight bounds on the worst
case complexity. More precisely, 
$$\text{comp}(\ve)\,=\,cam^{\KR}(\ve,\phi^{\mr}),\tag 5.2$$
where $c$ is the cost of one matrix-vector multiplication 
and 
$$a\in [0.5-1/m^{\KR}(\ve,\phi^{\mr})\,,1+10\,n/c].$$
For small $\ve$ and $c\gg n$, we have roughly $a\in 
[\tfrac 12,1]$. 

Because of (5.2), 
the problem of obtaining the complexity reduces to the 
problem of finding
$m^{\KR}(\ve,\phi^{\mr})$. 
This number is known for some classes 
$F$; see TW [84] and TWW [88, Chapter 5, \S 9]. 
We discuss two classes: 
$$\align F_1&=\big\{A:\ A=A^{\roman T}>0,\ 
\text{and}\ \|A\|_2\,\|A^{-1}\|_2\le M
\big\},\\ F_2&=\big\{A:\ A=A^{\roman T},\,\ \ \ \ \ \ 
\text{and}\ \|A\|_2\,\|A^{-1}\|_2\le M\big\}.\endalign$$
That is, $F_1$ is the class of symmetric positive 
definite matrices 
with condition numbers bounded uniformly by $M$. Here $M$ 
is a given number, 
$M\ge 1$. The class $F_2$ differs from $F_1$ by the lack 
of 
positive definiteness.

For these two classes, the result of Nemirovsky [91] can 
be applied and 
for $m(\ve)\le \tfrac 12(n-3)$ we have 
better bounds on $a$; namely 
$a\in [1-1/m^{\KR}(\ve,\phi^{\mr})\,,1+10\,n/c].$
Thus, for small $\ve$ and $c\gg n$, $a\simeq 1$. 

For the class $F_1$, we have 
$$m^{\KR}(\ve,\phi^{\mr})
=\min\left\{ n,\left \lceil \frac{
\ln\big((1+(1-\ve ^2)^{1/2})/\ve\big)}
{\ln\big((M^{1/2}+1)/(M^{1/2}-1)\big)} \right 
\rceil\right\} .$$
For small $\ve$, large $M$, and $n>M^{1/2}\ln\,(2/\ve 
)/2$, we have
$$m^{\KR}(\ve,\phi^{\mr})\simeq \frac{\sqrt{M}}2\ 
\ln\frac 2{\ve}.$$
For the class $F_2$, we have 
$$
m^{\KR}(\ve,\phi^{\mr})=\min \left\{n,2\left \lceil \frac 
{\ln ((1+(1-\ve ^2)^{1/2})/\ve)}
{\ln((M+1)/(M-1))} \right \rceil\right\}.
$$
For small $\ve$, large $M$, and $n>M\ln\,(2/\ve)$, we have
$$m^{\KR}(\ve,\phi^{\mr})\simeq M\ln\frac 2{\ve}.$$

These formulas enable us to compare the complexities for 
classes $F_1$
and $F_2$. For small~$\ve$, large $M$, and 
$n>2M\ln\,(2/\ve)+3$, we have 
$$
\frac {\text{comp}(\ve,F_1)}{\text{comp}(\ve,F_2)}\ \simeq
\ \frac 1{2\sqrt M}.$$ 
This shows how positive definiteness decreases the 
$\ve $-complexity.

P [92] has four ``high level'' criticisms of IBC research 
on the large linear systems problem. We also select three 
additional
criticisms from P [92, 4]. We shall respond to these 
seven criticisms. P [92] contains other misunderstandings 
and errors regarding this topic but we will not try the 
reader's patience by responding to each of these. 
We list the seven criticisms of P [92]:

\roster
\item"{(i)}" IBC ``is not complexity theory'' since ``the 
stubborn fact remains that restricting information to 
Krylov 
information is not part of the linear equations problem'' 
P~[92, 2.A].
\item"{(ii)}" ``The trouble with this apparent novelty is 
that it is
not possible to evaluate the residual norm $\|Az-b\|$ for 
those
external $z$ because there is no known matrix $A$ (only 
Krylov
information). So how can an algorithm that produces $z$ 
verify whether
or not it has achieved its goal of making 
$\|Az-b\|<\ve\|b\|$'' P
[92, 2.C].
\item"{(iii)}" ``The ingredient of IBCT that allows it to 
generate 
irrelevant results is the problem class $F$ [see 
paragraph 2 in (A)]. $F$
did not appear in our brief description of the theory in 
the second
paragraph of \S1 because it is not a logically essential
ingredient but rather a parameter within IBCT;'' P [92, 
2.B].
\item"{(iv)}" ``IBCT's suggestion that it goes beyond the 
well-known 
polynomial class of algorithms is more apparent than 
real;'' P [92,
2.C].
\item"{(v)}" ``Here is a result of ours that shows why 
the nonpolynomial 
algorithms are of no interest in worst case complexity;'' 
P[92, 4.3].
\item"{(vi)}" ``With a realistic class such as SPD (sym, 
pos. def.) MR
is optimal (strongly) as it was designed to be, and as is
well known;'' P [92, 4.4].
\item"{(vii)}" ``The theory claims to compare algorithms 
restricted
solely to information $N_j$. So how could the Cheb 
algorithm obtain
the crucial parameter $\rho$?;'' P [92, 4.4].
\endroster
\indent We respond to each of these seven  criticisms.

(i) IBC does {\it not } restrict information to Krylov 
information.
The optimality of Krylov information in the class of 
matrix-vector multiplications is a conclusion, not an 
assumption. 

IBC does assume a class $\Lambda$ of information 
operations. The reasons
why this is both necessary and beneficial were discussed 
in \S2.
Here we confine ourselves to certain classes relevant to 
large linear
systems. 

Let $\Lambda_1$ denote the class of matrix-vector 
multiplications. 
Then as described above, 
for an orthogonally invariant class $F$ we may {\it 
conclude} that 
Krylov information is optimal to within 
a multiplicative factor of at most 2. 
Furthemore, we may {\it conclude} that Krylov information 
and the mr 
algorithm are almost $\ve$-complexity optimal. 
Rather tight bounds have been obtained 
on the complexity of important classes such as $F_1$ and 
$F_2$, see 
above. Additional classes of matrices are studied in TW 
[84]. 

Let $\Lambda_2$ denote the class of information 
operations where inner 
products of rows (or columns) of $A$ and an arbitrary 
vector $z$ can be 
computed. Rabin [72] studied the class $\Lambda_2$ for 
the exact solution
of linear systems, $\ve=0$, and for an arbitrary 
nonsingular matrix~$A$.
He proved that, roughly, $\tfrac 12 n^2$ inner products 
are sufficient to solve
the problem. No results are known for $\ve>0$. 

Let $\Lambda_3$ denote the class of information 
operations consisting
of arbitrary linear functionals. Optimality questions for 
the class
$\Lambda_3$ are posed in TW [84].  No results are known and
we believe this to be a difficult problem. 

Let $\Lambda_4$ denote the class of information operations 
consisting of {\it continuous} nonlinear functionals, and 
let 
$\Lambda_5$ denote the class of nonlinear functionals. In 
general,
complexity results in $\Lambda_4$ and $\Lambda_5$ can be 
different;
see Kacewicz and Wasilkowski [86] and Math\'e [90]. For 
linear
systems, these classes are too powerful since 
all entries of the matrix $A$ can be recovered by knowing 
the value 
of one continuous nonlinear functional. Thus, the 
$\ve$-cardinality
number is $1$ even for $\ve=0$; see TW [80, Chapter 7, \S 
3] for 
related material. 

(ii) If the class $\Lambda$ consists of matrix-vector 
multiplications
then, of course, we can evaluate the residual $\|Az-b\|$ 
for any $z$. If 
$z$ is outside of a Krylov subspace this requires one 
additional 
matrix-vector multiplication. 

On the other hand, it is sometimes possible to guarantee 
that 
$\|Az-b\|\le \ve$, without 
computing the residual $\|Az-b\|$. This can be done 
by using a priori information that $A\in F$ and the 
computed 
Krylov information. An example of such a situation is 
provided by
the Chebyshev algorithm for the class 
$F=\{A=I-B:\ B=B^{\roman T},\,\|B\|\le\rho<1\}.$

In general, if the assumptions are satisfied, 
IBC is {\it predictive}. The results of the theory {\it 
guarantee} an
$\ve$-approximation. One simply does the amount of work 
specified by
the upper bound on the complexity. For important classes 
of matrices
we have seen above that there are rather tight bounds on 
the complexity. 
Therefore this strategy does not require much more work 
than necessary. 

For most problems there is {\it no} residual that can be 
checked. 
There are residuals for problems related to solving 
linear or nonlinear 
equations. In the multivariate integration example of \S 2,
there is no residual that can be computed. Yet, IBC 
guarantees an 
$\ve$-approximation by using a priori information about 
the class $F$. 

(iii) We responded in general to the criticism that $F$ 
is not 
needed in \S 4; here we focus on large linear systems. On 
this
problem P [92, 2.B]  states that ``IBCT seems to use $F$ 
as a 
tuning parameter designed to keep $k<n$.''

The domain $F$ is {\it not} a tuning parameter; it is 
needed
for the problem to be well defined. 
The domain $F$ contains {\it all} a priori knowledge 
about matrices $A$. 
The more we know a priori, the smaller the domain $F$ 
becomes, and
as $F$ becomes smaller, the problem becomes easier. 
Furthermore, a priori
information is often available in practice. For example, 
matrices
which occur in the approximation of elliptic partial 
differential
operators are symmetric positive definite, often  with 
known bounds
on condition numbers.

Fortunately, many important classes which occur in 
practice are 
orthogonally invariant and the $\ve$-complexity 
optimality of Krylov information and the mr algorithm may 
be applied.

Of course, numerical analysts use different algorithms 
for different
classes of matrices (symmetric, positive definite, 
tridiagonal,
Toeplitz, etc.)  It is therefore all the more surprising 
that P [92]
objects to the concept of the class $F$.

(iv) P [92, 2C] claims that there is no need to go 
``beyond the 
well-known polynomial class of algorithms.''
It should be obvious 
that {\it all} algorithms must be allowed to compete 
if we want to establish lower bounds on complexity. 

For orthogonally invariant classes it turns out that the 
restriction to 
the polynomial class of algorithms does not cause any harm 
since the classical mr algorithm may lose at most one
insignificant step. But this had to be proven!

In fact, it is not uncommon in computational complexity 
that the 
known algorithms (that use the specific information) 
turn out to be optimal or close to optimal.
Examples include the Horner algorithm for evaluating a 
polynomial, 
the finite element method with appropriate parameters for 
elliptic
partial differential equations, 
or the bisection algorithm for approximating a zero of a 
continuous
function that changes sign at the interval endpoints. 

For large linear systems, a sufficient condition 
for almost $\ve$-complexity 
optimality of Krylov information and 
the mr algorithm is orthogonal
invariance of the class $F$. As mentioned above, Example 
3.5 of
TW [84] shows that if $F$ is not 
orthogonally invariant, the mr algorithm may lose its 
optimality. 
In this example the restriction to the polynomial class 
of algorithms 
is harmful because the optimal algorithm is nonpolynomial. 

 (v) P [92, 4.3] supports his claim that nonpolynomial 
algorithms
are not interesting by the Theorem of \S 4.3. 
This theorem holds for the class of SPD of all $n\times 
n$ symmetric 
positive definite matrices. 
In this theorem it is shown that for every nonpolynomial
algorithm which computes 
an approximation outside the Krylov
subspace for $A \in \text{SPD}$, there exists 
a matrix from SPD which has the identical Krylov 
information as
$A$ and for which the residual is arbitrarily large.

We do not understand why the Theorem of \S 4.3 and
the one page sketch of its proof were supplied. The same 
statement can be 
found in Example 3.4 of TW [84]. In addition, Example 3.4 
shows that 
polynomial algorithms are also not good for the class 
SPD; that is,
$n$ matrix-vector multiplications are needed to compute an 
$\ve$-approximation. The reason neither polynomial nor 
nonpolynomial 
algorithms are good is that the class SPD is too large. 

We stress that Example 3.4 and the Theorem of \S 4.3 hold 
for
$F=$SPD. As mentioned above, for any orthogonally 
invariant class $F$ 
the nonpolynomial algorithms are not of interest since it 
has been 
proven that the mr algorithm is optimal, possibly modulo 
one 
matrix-vector multiplication. Also, as mentioned above, 
if $F$ is not 
orthogonally invariant, a nonpolynomial algorithm may be 
optimal. 

(vi) P [92] claims that the mr algorithm is optimal 
``as it was designed to be'' for the class SPD. 
This is simply not true. The mr algorithm is {\it 
defined} to be
optimal in the class of polynomial algorithms. Optimality
of the mr algorithm  in the class of {\it all} algorithms 
for the class SPD requires a proof. 

(vii) 
As already explained, the information that 
$A\in F=\{A=I-B:\ B=B^{\roman T},\,\|B\|\le\rho<1\}$ is 
{\it not } used by the mr algorithm. This means that the 
mr algorithm does not use the parameter $\rho$ which is 
assumed known {a 
priori} and may be used by competing
algorithms. The parameter $\rho$ is used by 
the Chebyshev algorithm and that is why the mr algorithm 
loses 
one step for the class $F$. P [92, 4.4] 
turns the positive
optimality result for the mr algorithm into the 
irrelevant question 
``how could the Chebyshev algorithm obtain the crucial 
parameter $\rho$?''
By the way, the parameter $\rho$ is {\it not} so crucial 
if it decreases 
the number of steps by {\it only} one!
\heading 6. Large eigenvalue problem
\endheading
P [92] has three ``high level'' criticisms of the IBC 
research 
on the large eigenpair problem. He also criticizes the 
numerical
testing. We shall respond to these four criticisms. 

We list the four criticisms of P [92]:
\roster
\item"{(i)}" Kuczy\'nski [86] computes an unspecified 
eigenvalue; 
P [92, 2.D].
\item"{(ii)}" IBC ``is not complexity theory.'' The 
reason given is
that ``the stubborn fact remains that restricting 
information to
Krylov information is not part $\ldots$ of the eigenvalue 
problem;''
P~[92, 2.A].
\item"{(iii)}" ``The fact that $b$ is treated as 
prescribed data is 
quite difficult to spot;'' P [92, 2.E].
\item"{(iv)}" ``The author has worked exclusively with 
tridiagonal 
matrices and has forgotten that the goal of the Lanczos 
recurrence is
to produce a tridiagonal matrix! Given such a matrix one 
has {\it no
need of either Lanczos or GMR};'' P [92, 5.5].
\endroster

We respond to each of these four criticisms.

(i) P [92] is certainly correct in asserting that when 
only 
one or a few eigenvalues of a symmetric matrix are 
sought, then one typically 
desires a preassigned eigenvalue or a few preassigned 
eigenvalues. To be 
specific, assume that the largest eigenvalue is to be 
approximated. 

It would be desirable to always guarantee that the 
largest eigenvalue 
$\lambda_1(A)$ of a 
large symmetric matrix $A$ can be computed to within 
error $\ve$. 
Unfortunately, this cannot be done with less than $n$ 
matrix-vector multiplications, that is, without 
recovering the matrix $A$; 
see TWW [88, Chapter 5, \S 10]. 
More precisely, let $F$ denote the class
of all $n\times n$ symmetric matrices and let $\Lambda$ 
consist of 
matrix-vector multiplications. That is, 
 $N(A)=[Az_1,\dots,Az_k]$,
where $z_1$ is an arbitrary vector and $z_i$ for $i\ge2$ 
may depend arbitrarily 
on $Az_1,\dots,Az_{i-1}$. Then for $k\le n-1$, there 
exists no 
such $N$ and no algorithm $\phi$ which uses $N$ such that 
$U(A)=\phi(N(A))$ 
satisfies
$$
|\lambda_1(A)-U(A)|\le \ve \|A\|,\quad \forall A\in F.
$$
We are surprised that although TWW [88] is cited in P 
[92], 
he does not seem to be aware of this result. 

Thus, the goal of computing an $\ve$-approximation to 
the largest eigenvalue of a large symmetric matrix cannot 
be achieved,
if less than $n$ matrix-vector multiplications are used. 
This is,
of course, a worst case result. There are a number of 
options for coping with 
this negative result. One could stay with the worst case 
setting but settle 
for an unspecified eigenvalue. Or one could give up on 
the worst case guarantee
and settle for a weaker one. We consider these options in 
turn. 

(i.1)
One option is to settle for an {\it unspecified 
eigenvalue}. 
More precisely, the problem 
studied by Kuczy\'nski [86] and Chou [87] is defined as 
follows. For $A\in F$,
compute $(x,\lambda)$ with $x\in \reals^n,\,\|x\|=1,$ and 
$\lambda \in \reals$,
such that 
$$
\|Ax-\lambda x\|\le \ve \|A\|.
$$

Chou proved, modulo a multiplicative factor of 2, 
optimality of Krylov
information $N(A)=[Ab,\dots,A^kb]$, where 
$b$ is a nonzero vector. Optimality of Krylov information 
holds 
independently of the choice of the vector $b$. 
Kuczy\'nski proved,
modulo an additive term of 2, optimality of the 
generalized minimal 
residual (gmr) algorithm that uses Krylov information. 
(Optimality of Krylov 
information and the gmr algorithm is understood as in \S5.
These optimality results hold for any orthogonally 
invariant class 
of matrices. )

Since the gmr algorithm has small combinatory cost, we 
conclude that 
the pair Krylov information and gmr 
algorithm is (almost) $\ve$-complexity optimal. 
Kuczy\'nski found good bounds on the worst case error of 
the gmr algorithm.
Hence, for $n>\ve^{-1}$, 
the worst case $\ve$-complexity is given by 
$$\text{comp}(\ve)=\frac {ac}\ve,
$$
where $a$ roughly belongs to $[\tfrac 14,1]$ and, as 
before, 
$c$ is the cost of one matrix-vector multiplication. 

(i.2)
A second option is to attempt to approximate the largest 
eigenvalue but to settle for a weaker guarantee.
KW [89]
\footnote{This paper is mistakenly referred in P [92] as 
[Tr $\&$ Wo, 1990].} 
study this problem in the randomized setting. (See, e.g., 
TWW [88, Chapter 11] for a general discussion of the 
randomized setting.)

In particular, the Lanczos algorithm is studied. The 
Lanczos 
algorithm uses Krylov information 
$N(A)=[Ab,A^2b,\dots,A^kb]$ with a 
{\it random} vector $b$ which is uniformly distributed 
over the unit sphere
of $\reals^n$. The error is defined  for a fixed matrix 
$A$ 
while taking the {\it average} with respect to the 
vectors $b$. 

To date only  an upper bound on the error of 
the Lanczos algorithm  with randomized Krylov information 
has been obtained. 
This upper bound  is proportional to  $((\ln n)/k)^2$. 

As always, to obtain complexity results both the 
information and 
the algorithm must be varied. Lower bounds are of 
particular interest. 
The complexity of approximating the largest eigenvalue 
in the randomized setting is open.

(ii) P [92, 2.A] states
``$\ldots$ the stubborn fact remains that restricting 
information to Krylov
information is not part $\ldots$ of the eigenvalue 
problem.''

Although we have mentioned several times in this paper 
that P [92]
seems unaware of the results regarding optimality of Krylov
information we are particularly surprised that he appears 
unaware of
this result in the context of the large eigenvalue 
problem.  P [92]
repeatedly cites Kuczy\'nski [86] where Chou's result is 
reported.

(iii) P [92, 2.E] states ``the fact that $b$ is treated 
as prescribed data is quite difficult to spot.''
Perhaps the reason it is difficult to spot is that it is 
not prescribed.

What is assumed known? It is known a priori that $A$ is a 
symmetric 
$n\times n$ matrix. Furthermore, we are permitted to 
compute $Az_1,\dots
,Az_k$, where $z_i$ may be adaptively chosen. We are 
permitted to choose
$z_1$, which is called $b$, {\it arbitrarily}. In 
choosing $b$ we cannot 
assume that $A$ is known, since the raison d'etre of 
methods for solving 
large eigenvalue problems is just that $A$ need not be 
known. 

By the result quoted in (i),  
it is impossible to guarantee that we can find a vector 
$b$ such that an 
$\ve$-approximation to the largest eigenvalue can be 
computed for
all symmetric $n\times n$ matrices with $k<n$. 

If Krylov information $Ab,A^2b,\dots,A^kb$ is used then 
the situation is 
even worse. Even for arbitrary $k$, i.e., even for $k\ge 
n$, 
an $\ve$-approximation cannot be computed. 
Indeed, suppose we choose a 
vector $b$ and a matrix $A$ such that  $Ab=b$. 
Then Krylov information is reduced just 
to the vector $b$. The largest eigenvalue cannot be 
recovered (unless $n=1$). 
Thus, for any vector $b$ there are symmetric 
matrices $A$ for which Krylov information will not work.

Of course, one can choose $b$ randomly, as was discussed 
above.
The average behavior with respect to vectors $b$ is 
satisfactory 
for {\it all} symmetric matrices. 
But then one is settling for a weaker guarantee of 
solving the problem. 

P [92, 2.E] claims that for Krylov information 
``satisfactory starting
vectors are easy to obtain.''  This remark seems to 
confuse the worst
case and randomized settings. To get a satisfactory 
starting vector
$b$ in the worst case setting, the vector $b$ must be 
chosen using
some additional information about the matrix $A$. If such 
information
is not available, it is impossible to guarantee 
satisfactory starting
vectors.  On the other hand, in the randomized setting it 
is indeed
easy to get satisfactory starting vectors.

(iv) P [92, 5.5] complains that Kuczy\'nski [86] 
tests only tridiagonal matrices. 

There is no loss of generality in  restricting 
the convergence tests of the Lanczos or gmr algorithms to 
tridiagonal 
matrices. That was done in Kuczy\'nski [86] to speed up 
his tests. 
What is claimed 
in  Kuczy\'nski [86] for the pairs $(TRI,b)$, $TRI$ a 
tridiagonal matrix and 
$b=e_1=[1,0,\dots,0]^{\roman T}$,
is also true for the pairs $(Q^{\roman 
T}\,TRI\,Q,Q^{\roman T}b)$ 
for any orthogonal matrix $Q$.
Obviously, the matrix $Q^{\roman T}\,TRI\,Q$ is {\it 
not}, in general, tridiagonal.

The confusion between the worst case and randomized 
settings is also 
apparent when P~[92] discusses numerical tests performed 
by 
Kuczy\'nski [86]  and by him.

For the unspecified eigenvalue problem, 
Kuczy\'nski [86] compares the gmr and Lanczos algorithms 
in the 
{\it worst case setting}. These two algorithms cost 
essentially the same
per step, and the gmr algorithm never requires more steps 
than the
Lanczos algorithm. For some matrices, 
the gmr algorithm uses substantially fewer steps than the 
Lanczos 
algorithm. That is why in the {\it worst case setting} 
the gmr algorithm
is preferable.

P [92] performed his numerical tests for the Lanczos 
algorithm with
{\it random} starting vectors $b$. Thus, he uses a 
different setting.
It is meaningless to compare numerical results in 
different settings.

Finally, extensive numerical testing is also
reported in KW [89] for 
approximating the largest eigenvalue by the Lanczos 
 algorithm 
with randomized starting vectors. The Lanczos algorithm 
worked quite well
for all matrices tested. The numerical tests
reported by P [92] and KW [89] show the efficiency of
the Lanczos algorithm in the randomized setting. 

\heading  7. Refinements of IBC
\endheading
Our response to the criticism in P [92] does not mean 
that the
current model assumptions of IBC are the only ones 
possible. On the 
contrary, we believe that in some circumstances these 
assumptions 
should be refined to improve the modelling of 
computational problems. 
We have mentioned the desirability of such refinements 
in, e.g., 
TWW [88, Chapter 3, \S2.3]
and W [85, \S 9]. 
In this section we will very briefly indicate some of the 
possible refinements and extensions of IBC, and indicate 
partial
progress. This is preparatory to responding to several 
comments in P [92]. 

Refinements and extensions of IBC include  the following:
\roster
\item We usually assume the real number model in a 
sequential model of
computation where the cost of a
combinatory operation is independent of the precision of 
the operands
or of the result. Also of interest is a model where the 
cost of a 
combinatory operation depends on the precision (bit 
model) and/\!or 
on the particular operation. Parallel and distributed 
models of 
computation should also be studied. For examples of work 
in these
directions see
Boja\'nczyk [84] who studies the approximate solution of 
linear
systems using a variable precision parallel model of 
computation,  
and Kacewicz [90] who studies initial value problems for 
both
sequential and  parallel models of computation.
\item  We usually assume that for every information 
operation $L\in 
\Lambda$ and for every $f\in F$ the computation of 
$L(f\,)$ costs $c$, $c>0$.
Also of interest is a model where the cost of an 
information operation
depends on $L$, $f$, and precision. For an example, see
Kacewicz and Plaskota [90] who study 
linear problems in a model where the cost of information 
operations varies with the computed precision. 
\item Let $S$ be a linear operator. Then we often assume 
that the set
$F$ is balanced and convex; TWW [88, Chapter 4, \S 5]. In
particular, for functions spaces, we often assume that 
$F$ is a
Sobolev space of smoothness $r$ with a uniform bound on 
$\|f^{\,(r)}\|$. It is of interest to study $F$ which do 
not have 
such a nice structure.
\endroster

P [92, 1] states ``a handful of reservations about IBC
have appeared in print.'' These ``reservations'' turn out 
to concern
refinements of IBC.  P [92] writes that Babu\u ska [87] 
calls for
realistic models. For example, Babu\u ska points out that 
for some
problems arising in practice the set $F$ does not consist 
of smooth
functions but rather of functions which are piecewise 
smooth with 
singularities at unknown points. We agree that this is an 
important
problem. A promising start has been made by 
Wasilkowski and Gao [89] on estimating a singularity of  
a piecewise smooth function in a probabilistic setting. 

Babu\u ska observes that the user may not know the class 
$F$ or not
know $F$ exactly, and suggests the importance of algorithms
which enjoy optimality properties for a number of 
classes.  We agree 
that this is an important concern and a good direction 
for future
research. See W [85, \S 9.3] where this problem is called 
the ``fat'' $F$ problem and where 
partial results are discussed. One attack on
this problem is to address 
the domain membership problem defined in \S 4. As indicated
there, this can only be done with a stochastic assurance. 

P [92, 1] asserts that in a review of TWW [83], Shub [87] 
 ``gives a couple 
of instances of unnatural measures of cost.'' (These 
 words are from P
[92], not from Shub [87].) Shub, 
in a generally favorable review (the reader may want to 
verify this), suggests circumstances when the cost of an 
information operation should vary. We concur.

\heading 8. Summary
\endheading
P [92, 2] states five high level criticisms of IBC. We 
responded
to them in the following sections:

$$\matrix \format \c&\quad\c\\
\text{Criticism}&\text{Response}\\
\text{A}&\text{1,\ 3,\ 5,\ 6}\\
\text{B}&\text{4,\ 5,\ 6}\\
\text{C}&\text{5}\\
\text{D}&\text{6}\\
\text{E}&\text{6}
\endmatrix
$$	
There are additional criticisms, and in \S\S 5 and 6 
we responded to the ones which seem most important.

P [92, 1] states that ``a handful of reservations about 
IBCT have
appeared in print.'' He neglects mentioning the many 
favorable
reviews.  He cites two examples of reservations.  We 
discussed the
comments of Babu\u ska [87] and Shub [87] in \S 7.

P [92] is based upon the following syllogism:
\roster
\item {\bf Major Premise:} If two specific papers of IBC 
are
misleading, then IBC is flawed.
\item {\bf Minor Premise:} Two specific papers of IBC 
regarding  matrix computations are misleading.
\item {\bf Conclusion:   } IBC is flawed. 
\endroster

We have shown that his reasons for believing the minor 
premise
are mistaken.

\enddocument